\documentclass[a4paper]{amsart}

\usepackage[latin1]{inputenc}
\usepackage[T1]{fontenc}
\usepackage{setspace}
\usepackage{color}
\usepackage{amsmath, graphicx, cite, enumerate, comment}
\usepackage{amssymb, amsthm, amscd}
\usepackage{latexsym}
\usepackage[T1]{fontenc}	
\usepackage[english]{babel}
\usepackage{amsfonts}
\usepackage{amscd}
\usepackage{cite}

\newtheorem{theorem}{Theorem}[section]

\newtheorem{corollary}[theorem]{Corollary}

\newtheorem{proposition}[theorem]{Proposition}

\theoremstyle{definition}
\newtheorem{definition}[theorem]{Definition}

\theoremstyle{remark}
\newtheorem{remark}[theorem]{Remark}

\numberwithin{equation}{section}

\begin{document}
\date{}
\title[Finiteness results for minimal surfaces]
{Generic finiteness of minimal surfaces with bounded Morse index}
\author{Alessandro Carlotto}
\address{Institute for Theoretical Studies \\
	ETH \\
	Z\"urich}
\email{alessandro.carlotto@eth-its.ethz.ch}

\begin{abstract}
Given a compact 3-manifold $N$ without boundary, we prove that for a bumpy metric of positive scalar curvature the space of minimal surfaces having a uniform upper bound on the Morse index is always finite unless the manifold itself contains an embedded minimal $\mathbb{R}\mathbb{P}^2$. In particular, we derive a generic finiteness result whenever $N$ does not contain a copy of $\mathbb{R}\mathbb{P}^3$ in its prime decomposition. We discuss the obstructions to any further generalization of such a result. When the metric $g$ is required to be (scalar positive and) strongly bumpy (meaning that all closed, \textsl{immersed} minimal surfaces do not have Jacobi fields, a notion recently proved to be generic by B. White) the same conclusion holds true for \textsl{any} closed 3-manifold.
\end{abstract}	

\maketitle

\maketitle

\section{Introduction} \label{sec:intro}



Let $(N,g)$ be a compact, boundaryless Riemannian manifold of dimension three: then we know that under various sorts of natural topological or curvature  assumptions $N$ contains infinitely many distinct closed, embedded minimal surfaces. For instance, this is definitely the case for certain classes of aspherical 3-manifolds whose fundamental group contains a (finitely generated) noncyclic abelian subgroup (by virtue of the work of Schoen-Yau \cite{SY79} and Freedman-Hass-Scott \cite{FHS83}, see also Meeks-Simon-Yau \cite{MSY82}) or on the other hand if $Ric_g>0$ (by virtue of the work of Marques-Neves, see \cite{MN13}). However, it is natural to wonder whether the subclass of closed, embedded minimal surfaces of \textsl{bounded complexity} is compact (with respect to single-sheeted graphical convergence in $C^{k}$ for $k\geq 2$) or perhaps even finite. Here the phrase \textsl{bounded complexity} may be specified to assume (at least) two distinct, and a priori inequivalent meanings. These correspond to the following sublcasses:
\begin{enumerate}
\item{$\mathcal{G}(c)=\left\{ M\subset N \ \textrm{minimal surface with} \ genus(M)\leq c \right\}$}
\item{$\mathcal{I}(I)=\left\{ M\subset N \ \textrm{minimal surface with} \ Ind(M)\leq I \right\}$}
\end{enumerate}	
where (in both cases) $M$ is assumed to be closed and embedded and $Ind(M)$ is the Morse index of $M$, namely the index of the stability operator of $M$ as a minimal surface in $(N,g)$.

On the compactness side, it has been known for about three decades, namely since the work by Choi and Schoen \cite{CS85}, that the space of closed, embedded minimal surfaces of fixed topological type in $N$ is compact in the $C^{k}$ topology for any $k\geq 2$ provided the Ricci curvature of $N$ is positive. In fact, combining this result with those in \cite{Sha15} and in \cite{ACS15} one can derive the following rather complete description of the scenario:

\begin{theorem}(Corollary 1.7 in  \cite{ACS15})\label{cor:3dscenario}
	Let $\mathfrak{C}\subset\mathfrak{M}$ be a subclass of closed minimal surfaces inside some smooth closed Riemannian manifold $(N,g)$ of dimension $3$ satisfying $Ric_g >0$. Then a uniform bound on any one of the following quantities for every $M\in\mathfrak{C}$ is enough to ensure compactness of $\mathfrak{C}$ and hence leads to a bound on the rest of them for every $M\in\mathfrak{C}$:
	\begin{itemize}
		\item the genus of $M$ 
		\item $Ind(M) + \mathcal{H}^2(M)$
		\item $\lambda_p(M) + \mathcal{H}^2(M)$
		\item $\sup_M |A| + \mathcal{H}^2(M)$
		\item $\int_M |A|^2+\mathcal{H}^2(M)$. 
	\end{itemize}
Here $A$ is the second fundamental from, while $\lambda_p(M)$ denotes the $p-$th eigenvalue of the stability operator of $M$.	
\end{theorem}

We shall refer the reader to the Introduction of \cite{ACS15} for strong compactness results when the positivity assumption on the ambient Ricci curvature is relaxed or even dropped.

 Unfortunately, very elementary examples show that finiteness cannot be expected in general to hold for any of the classes introduced above: for instance, if the isometry group of $(N,g)$ is transitive then one can choose the parameters $c$ and $I$ so that the corresponding classes $\mathcal{G}(c),\mathcal{I}(I)$ are not empty and hence have infinite cardinality.

A more reasonable question is then whether finiteness of those classes can be \textsl{generically} true, namely for a generic choice of the ambient Riemannian metric $g$. In this work, we shall deal with the following natural notion of genericity: given a Riemannian 3-manifold $(N,g)$ as above we shall say that the metric $g$ is bumpy, if any closed \underline{embedded} minimal surface has no non-trivial Jacobi fields. A well-known application of Sard's Lemma in a Banach space setting, due to B. White \cite{Whi91}, ensures that the family of bumpy metrics on a given smooth 3-manifold is residual (in the sense of Baire category). 

In this article, we provide a criterion that ensures generic finiteness of the class $\mathcal{I}(I)$ for any choice of $I\geq 0$. More precisely, our main result is the following statement.

\begin{theorem}\label{thm:main}
Let $N$ be a compact, orientable 3-manifold\footnote{Throughout this paper all 3-manifolds are assumed to be orientable.} without boundary and let it be endowed with a bumpy Riemannian metric $g$ of positive scalar curvature. Then the following dichotomy holds: \underline{either} for any integer $I\geq 0$ the set $\mathcal{I}(I)$ is finite 		
\underline{or} $N$ contains a minimally embedded copy of $\mathbb{R}\mathbb{P}^2$.
\end{theorem}	

Now, we know by the blow-up analysis performed by Perelman on Hamilton's Ricci flow \cite{Per02, Per03a, Per03b} (but see also earlier work by Gromov-Lawson \cite{GL80} and Schoen-Yau \cite{SY79}) that a Riemannian manifold $(N,g)$ of positive scalar curvature has a prime decomposition where each summand is either a copy of $S^2\times S^1$ or a spherical space form. On the other hand, a simple topological argument shows that if a 3-manifold $(N,g)$ contains an embedded $\mathbb{R}\mathbb{P}^2$ (minimal or not) then one can find a regular neighborhood of such surface, say $U$, which is diffeomorphic to $\mathbb{R}\mathbb{P}^3$ minus a ball and hence, as an immediate consequence, $N$ must have a copy of $\mathbb{R}\mathbb{P}^3$ in its prime decomposition. Therefore, we can derive the following corollary. 

\begin{corollary}\label{cor:simple}
	Let $N$ be a compact, orientable 3-manifold without boundary not having any copy of $\mathbb{R}\mathbb{P}^3$ in its prime decomposition and let $g$ be a bumpy Riemannian metric of positive scalar curvature on $N$. Then for any integer $I\geq 0$ the set $\mathcal{I}(I)$ is finite.	
\end{corollary}	

Of course, the topological condition in question can be read off at the level of the fundamental group of $N$ and corresponds to the requirement that such group does not contain any subgroup of order two. 
We explicitly remark that Corollary \ref{cor:simple} immediately implies, as a byproduct, that the set $\mathcal{G}(c)$ is also generically finite for any $c\geq 0$ in the subclass of metrics of positive Ricci curvature (due to the fact that a bound on the genus does imply a bound on the Morse index thanks to the work of Eijiri-Micallef \cite{EM08} and Choi-Wang \cite{CW83}). 

\

Further extension of our theorem are limited by the following two classes of examples:
\begin{enumerate}
\item{if $N=S^1\times S^1\times S^1$ and $g$ is \textsl{any} metric, then $N$ contains infinitely many distinct, embedded minimal surfaces (each of them obtained by minimizing the area functional in a suitable homotopy class) and so finiteness for $\mathcal{I}(I)$ cannot possibly hold for all closed 3-manifolds under a sole scalar curvature bound of the form $R \geq \underline{R}, \ \underline{R}\in\mathbb{R}$;}	
\item{if $N$ is \textsl{any} three-manifold, then by virtue of \cite{CM99} we know the existence of an open set $\mathcal{O}$ of Riemannian metrics such that for any $g\in\mathcal{O}$ the manifold $(N,g)$ contains infinitely many distinct, embedded stable minimal tori, hence even if $N$ is simply-connected one cannot expect finiteness of $\mathcal{I}$ if the ambient metric is not scalar positive.}	
\end{enumerate}	
In particular, the enlightening examples provided by Colding and Minicozzi ensure that even when $N$ is simply-connected (namely $N=S^3$ by virtue of the positive solution of the Poincar\'e conjecture) there are choices of suitable negatively curved metrics that imply the existence of stable, minimal surfaces of fixed topological type, but arbitrarily large area.

We shall further remark here that according to Yau (see \cite{Yau82}) $S^3$ contains infinitely many closed minimal surfaces whatever Riemannian metric it is endowed with, thus our main theorem provides some information about the asymptotic behaviour of any sequence of such surfaces.

\begin{corollary}
(subject to Yau's conjecture)
Let $g$ be a bumpy Riemannian metric on $S^3$ of positive scalar curvature. Then there exist embedded minimal surfaces of arbitrarily large Morse index in $(S^3,g)$.	
\end{corollary}	 

This is somehow unexpected, as when the Ricci positivity is dropped Frankel's theorem is no longer true and hence the occurrence of finitely many models of min-max minimal surfaces with higher and higher multiplicity would not contradict the Gromov-Guth bounds for the $n-$width $\omega_n$ of $(S^3,g)$  with the consequence that min-max methods may in fact not be able to provide minimal surfaces of large Morse index (as in \cite{MN13}). 

\

These sorts of finiteness questions have attracted considerable interest and proved to be rather subtle, even in the most fundamental case concerning simple closed geodesics on closed surfaces for which we refer the reader to the beautiful papers by Colding and Hingston \cite{CH03, CH06} and reference therein.

Coming back to closed surfaces in 3-manifolds, our work was preceded by \cite{CM00} and, much more recently, by \cite{LZ15} where a generic finiteness result for $\mathcal{I}(I)$ in a Riemannian manifold $(N,g)$ of positive Ricci curvature was proposed. Our treatment, while identifying a topological obstruction, extends such result under a much less restrictive curvature condition on the ambient metric and, correspondently, requires some ideas and techniques which are not present in those works.
However, the logical structure of our proof goes back, in its essentials, to \cite{CM00} where the authors prove a finiteness result for the space 
\[
\mathcal{G}(\Lambda, c)=\left\{ M\subset N \ \textrm{minimal surface with} \ \mathcal{H}^2(M)\leq \Lambda, \ genus(M)\leq c \right\}.
\]
Suppose by contradiction the set in question had infinite cardinality: then it would contain a sequence $\left\{ M_i\right\}$ for which one can make use of the area bound and the genus bound (which, in turn, easily implies a uniform bound on the total curvature) to prove convergence to a limit minimal surface $M$, the convergence being smooth (possibly with multiplicity) away from finitely many exceptional points. At that stage, normalizing the height function of $M_i$ over $M$ (or between adjacent leaves of $M_i$) allows to construct a global solution of the Jacobi equation $Lu=0$, thereby contradicting the bumpyness assumption. In our case, we deal with two significant and strictly intertwined complications: first, we do not have area bounds, so that we cannot expect strong (possibly multiply-sheeted) convergence to hold but need to take a limit in the weaker sense of laminations\footnote{This delicate issue is also present in the work of Li-Zhou \cite{LZ15}, even though our treatment is partly different from theirs.} (as per \cite{CM04}); second, in this process of passing to the limit we may a priori generate minimal leaves that are complete but not compact, namely that spiral inside $N$ in a wild fashion. Our assumption on the positivity of the scalar curvature is a priori not incompatible with this occurrence (nor is of course the presence of closed stable minimal surfaces). At that crucial stage, our proof exploits an interesting $H-$radius estimate due to Schoen and Yau \cite{SY83}: we are deeply indebted to R. Schoen for pointing out the potential relevance of that result with respect to our work.  
In fact, we find that the statement of Proposition \ref{prop:Hradsurf} is of independent interest and its range of applicability could certainly go beyond the specific scopes of this paper.

\

Let us now comment on the relevance of the requirement that $N$ does not contain any minimally embedded $\mathbb{R}\mathbb{P}^2$. It is well-known that the lift of a stable $\mathbb{R}\mathbb{P}^2$ under the standard covering map $S^3\to\mathbb{R}\mathbb{P}^3$ is a minimal $S^2$ of Morse index 1, namely this negative direction ``disappears"  when passing to the quotient. Essentially, the same mechanism lies behind the assertion that the lift (under a Riemannian covering map $\pi: \tilde{N}\to N$) of a bumpy metric need not be bumpy.  Because of this fact, we are inclined to believe that there should exist bumpy Riemannian metrics on $\mathbb{R}\mathbb{P}^3$ (of positive scalar curvature) which contain infinitely many stable minimal spheres converging to a minimal $\mathbb{R}\mathbb{P}^2$ with two-sheeted graphical convergence. This phenomenon is described in more detail in Remark \ref{rem:counter} and seems to show that the statement of Theorem \ref{thm:main} is completely sharp and captures an essential topological restriction.
On the other hand, one could consider a slightly different notion of bumpyness: given a Riemannian 3-manifold $(N,g)$ as above we shall say that the metric $g$ is strongly bumpy, if any closed \underline{immersed} minimal surface has no non-trivial Jacobi fields. This notion has been proven to be generic (more precisely: Baire residual) in the very recent article \cite{Whi15} by B. White. Now, if a closed 3-manifold $(N,g)$ contains a minimally embedded $\mathbb{R}\mathbb{P}^2$ then there is of course a minimal (two to one) immersion $\varphi:S^2\to N$ and hence the argument we are about to present for Theorem \ref{thm:main} leads to the following conclusion. 

\begin{theorem}\label{thm:main2}
	Let $N$ be a compact, orientable 3-manifold without boundary and let it be endowed with a strongly bumpy Riemannian metric $g$ of positive scalar curvature. Then for any integer $I\geq 0$ the set $\mathcal{I}(I)$ is finite. 		
\end{theorem}	We shall remark here that this result has also been proven, independently of us and roughly at the same time (but with somewhat different arguments) by Chodosh-Ketover-Maximo \cite{CKM15}. In the context of strongly bumpy metrics, Theorem \ref{thm:main2} is patently sharp, as is witnessed by the aforementioned examples (1) and (2) above.

\

\textsl{Acknowledgments}. The author wishes to express his gratitude to Prof. Andr\'e Neves and to Prof. Richard Schoen for a number of enlightening conversations and would like to thank Prof. Davi Maximo for pointing out the article \cite{CM99}. The extensive feedback provided by Lucas Ambrozio and Ben Sharp on earlier versions of this paper has helped the author improve its quality in a dramatic fashion: he is deeply indebted to both of them. Thanks also to Prof. Richard Bamler and Prof. William Minicozzi for important clarifications and to the anonymous referee for accurate proofreading.
This work was done while the author was an ETH-ITS fellow: the outstanding support of Dr. Max R\"ossler, of the Walter Haefner Foundation and of the ETH Zurich Foundation are gratefully acknowledged.

\section{A collection of ancillary results}\label{sec:prelim}
	
\subsection{Jacobi operator and stability}	

We shall recall here the definition of the Morse index and the Jacobi eigenvalues $\lambda_k$ for general smooth minimal  surfaces $M\hookrightarrow (N,g)$. First of all, if $M$ is orientable then the second variation of the area functional can be written down purely in terms of sections of the normal bundle $v\in \Gamma(\textrm{Nor}(M))$ by
\[Q(v,v):=\int_{M} |\nabla^\bot v|^2 - |A|^2|v|^2 -Ric_g(v,v).
\]
Standard results on the spectra of compact self-adjoint operators on separable Hilbert spaces tell us that there is an orthonormal basis $\{v_i\}_{i=1}^{\infty}$ of $L^2(\Gamma(\textrm{Nor}(M)))$ consisting of eigenfunctions for the operator
\[
L^{\bot}v:=\Delta^{\bot} v +|A|^2 v + Ric^{\bot}_g(v)
\]
 with associated eigenvalues $\{\lambda_i\}_{i=1}^{\infty}$ of $Q$.

Now, if $M$ is non-orientable then we simply lift the problem to its orientable double cover $\tilde{M}$ via $\pi :\tilde{M}\to M$. Consider the linear subspace of smooth sections $v\in\Gamma(\pi^{\ast}(\textrm{Nor}(M)))$ such that $v\circ \tau =v$ where $\tau:\tilde{M}\to \tilde{M}$ is the unique deck transformation of $\pi$ which reverses orientation. Denote this subspace by $\tilde{\Gamma}(\pi^\ast (\textrm{Nor}(M)))$. We can also pull back the quantities $|A(x)|^2:= |A(\pi(x))|^2$ and $Ric_g(a,b) := Ric_g (\pi_\ast a, \pi_\ast b)$. Thus consider the quadratic form
\[\tilde{Q}(v,v):=\int_{\tilde{M}} |\nabla^\bot v|^2 - |A|^2|v|^2 -Ric_g(v,v) 
\]
over $\tilde{\Gamma}(\pi^\ast (\textrm{Nor}(M)))$. As before we can define the spectrum, and therefore index of $M$ to be that of $\tilde{M}$ with respect to $\tilde{Q}$ and $\tilde{\Gamma}(\pi^\ast (\textrm{Nor}(M)))$.

\

If the ambient manifold $N$ is orientable, a closed surface is orientable if and only if it is two-sided, namely if there exists a global section $\nu=\nu_M$ of its normal bundle inside  $\textrm{Tan}(N)$. If this is the case (namely if $M\hookrightarrow N$ is minimal and two-sided) the spectrum defined above patently coincides with the spectrum  of the \textsl{scalar} Jacobi or stability operator of $M$, namely 
\[
Lu:=\Delta_M u + (|A|^2+Ric_g(\nu,\nu))u
\]
 when we regard $L:W^{1,2}(M)\to W^{-1,2}(M)$.
\

\subsection{Minimal laminations and compactness}

In order to state the compactness theorem we are about to use in the sequel of this article, we shall first recall the definition of lamination, a notion that dates back to Thurston (see the most recent edition \cite{Thu97} of the well-known 1980 notes).

\begin{definition}
A codimension one \textsl{lamination} of (a three-dimensional manifold) $N$ is a collection $\mathcal{L}$ of smooth disjoint connected surfaces (called leaves) such that $\cup_{\Lambda\in\mathcal{L}}\Lambda$ is closed and, furthermore, for each $x\in N$ there exists an open neighborhood $U$ of $x$ and an associated local coordinate chart, $(U,\Phi)$, with $\Phi(U)\subset\mathbb{R}^3$ such that in these coordinates the leaves of $\mathcal{L}$ pass through the chart in slices of the form $(\mathbb{R}^2\times\left\{t \right\})\cap \Phi(U)$. 		
\end{definition}	

Therefore, we recall that a foliation is just a lamination without gaps, namely a lamination where the union of the leaves is all of $N$. We shall say that a lamination is minimal if each of its leaves is a minimal surface. Moreover, we shall agree that a sequence of laminations $\left\{\mathcal{L}_i\right\}$ \textsl{converges} to $\left\{\mathcal{L}\right\}$ if this is true, locally, for the corresponding coordinate charts (coherently, the class of convergence is determined by the functional spaces where the convergence of charts happens, in the obvious fashion).

In their monumental work on the Calabi-Yau conjectures, Colding and Minicozzi have developed powerful methods to take the limit of sequences of minimal surfaces for which no a priori area bounds are known: the limit objects that arise in such scenario are precisely minimal laminations (possibly with singularities, when no uniform curvature estimates hold). We refer the reader to the beautiful introductory article \cite{CM00} for an accessible presentation of some of these results and to the references therein for full details (for our purposes, the sole appendix B of the paper \cite{CM04} is the truly relevant part).

\begin{theorem}\label{thm:lamcomp}(essentially Prop. B.1 in \cite{CM04})
Let $(N,g)$ be a fixed 3-manifold. If $\left\{\mathcal{L}_i\right\}\subset B_{2r}(x)\subset N$ is a sequence of minimal laminations with uniformly bounded curvatures and each leaf has boundary contained in $\partial B_{2r}(x)$, then there exists a subsequence which converges in $B_{r}(x)$ in the $C^{0,\alpha}$ topology for any $\alpha$ to a Lipschitz lamination $\mathcal{L}$. The leaves of such lamination are smooth and minimal, and the convergence happens in $C^{\infty}$ leafwise.
\end{theorem}	

A few comments about this statement are appropriate. First of all, the fact that the convergence is smooth at the level of the single leaves is remarked (for instance) at page 115 of \cite{CM00}, in any case it is apparent from the proof given in \cite{CM04}. Indeed, the seemingly modest convergence (Lipschitz charts for $\mathcal{L}$ and $C^{0,\alpha}$ convergence) is only related to the \textsl{transversal directions} in view of the somehow restrictive definition of lamination, namely to the fact that one could have in $B_{r}(x)$ leaves of $\mathcal{L}_i$ which are graphs with (uniformly) small gradient, but tilted one with respect to the adjacent one(s). Hence, one might not expect too much control on a coordinate chart $\Phi:B_{r}(x)\to\mathbb{R}^3$ which is required to describe, \textsl{at the same time}, all leaves as horizontal planes. 

In order to handle potentially spiraling leaves it is useful to introduce the following definition.

\begin{definition}\label{def:conc}
Let $(N,g)$ be a Riemannian 3-manifold and let $\mathcal{L}$ be a minimal lamination. We say that $p\in\cup_{\Lambda\in\mathcal{L}}\Lambda$ is a concentration point for $\mathcal{L}$ if $\mathcal{L}$ has infinite density at $p$ (at all scales), namely if
\[
\liminf_{r\to 0}\frac{\mathcal{H}^2(B_r(p)\cap(\cup_{\Lambda\in\mathcal{L}}\Lambda))}{\pi r^2}=\infty.
\]   	
Else, if that is not the case, we say that 	$p\in\cup_{\Lambda\in\mathcal{L}}\Lambda$ is a non-concentration point for $\mathcal{L}$. 
\end{definition}

\begin{remark}\label{rem:equivconc}
Let $(N,g)$ be a Riemannian 3-manifold, let $\mathcal{L}$ be a minimal lamination and assume $\left\{\mathcal{L}_i\right\}$ is a sequence of minimal laminations converging to $\mathcal{L}$. Then it is immediately checked that $(1)\Rightarrow(2)$ where $(1)$ and $(2)$ are the following two assertions:
\begin{enumerate}
\item{$p\in\cup_{\Lambda\in\mathcal{L}}\Lambda$ is a concentration point for $\mathcal{L}$;}
\item{$p\in\cup_{\Lambda\in\mathcal{L}}\Lambda$ is an area-concentration point for the sequence $\left\{\mathcal{L}_i\right\}$, by which we mean that
\[
\liminf_{r\to 0}\liminf_{i\to\infty}\frac{\mathcal{H}^2(B_r(p)\cap(\cup_{\Lambda\in\mathcal{L}_i}\Lambda))}{\pi r^2}=\infty.
\]}	
\end{enumerate}		
\end{remark}	

\begin{remark}\label{rem:geomHarnack}
In the setting of Definition \ref{def:conc}, the geometric Harnack principle (see pp. 332-334 in \cite{Sim87}) implies that: if $\Lambda\in\mathcal{L}$ is a leaf (namely a connected component of the lamination in question) then \underline{either} all points of $\Lambda$ are concentration points for $\mathcal{L}$ (in which case we say that $\Lambda$ is an accumulating leaf) \underline{or} all points of $\Lambda$ are non-concentration points for $\mathcal{L}$ (in which case we say that $\Lambda$ is a non-accumulating leaf). 	
\end{remark}

\subsection{Removing isolated singularities}

In the proof of our main theorem, we shall construct a minimal lamination $\mathcal{L}$ by taking the limit of a sequence of closed minimal surfaces of bounded Morse index: as will be apparent from our argument the lack of uniform curvature estimates at finitely many points implies, \textsl{a priori}, the presence of punctures in $\mathcal{L}$ which we can get rid of by means of the following well-known removable singularity theorem.

\begin{theorem}\cite{GL86, CM04, MPR13}\label{thm:remsing}
Let $(N,g)$ be a complete 3-manifold and fix a point $p\in N$. Let $\mathcal{L}$ be a minimal lamination of $B_r(p)\setminus\left\{p\right\}$ for some $r>0$ and assume, furthermore, that each leaf of $\mathcal{L}$ has stable universal cover. Then $\mathcal{L}$ extends to a smooth minimal lamination of $B_r(p)$.
\end{theorem}	

We shall add here some words about the references above: while Theorem \ref{thm:remsing} would follow at once from Theorem 1.1 in \cite{MPR13} (because of the stability estimate by Schoen \cite{Sch83}), invoking the full strength of such result is not at all needed. Indeed, one can provide a proof relying on a blow-up argument and the classic Gulliver-Lawson Bersntein-type theorem \cite{GL86} (see also Lemma A.26 in \cite{CM15}). The key point is gaining properness (namely the fact that each connected component $\Lambda\subset B_r(p)\setminus\left\{p\right\}$ is properly embedded), for at that stage one can consider (by virtue of the monotonicity formula) the tangent cone of $M$ at $p$, prove that $M$ has finite Euler characteristic  and the conclusion is well-known (see e.g. Proposition 1 in \cite{CS85}). The argument which proves properness is just a variation of Proposition 4.2 in \cite{CCE15}, despite the the seemingly different scenario.  

\subsection{A double covering construction}

When the ambient manifold $(N,g)$ is not simply connected we may have to face the situation when some of the leaves of the minimal lamination $\mathcal{L}$ are not orientable, in which case the following covering space conctruction, due to X. Zhou, will be remarkably useful.

\begin{proposition} \label{pro:xz}(Proposition 3.7 in \cite{Zho12})
For any non-orientable embedded hypersurface $M^n$ in an orientable manifold $N^{n+1}$, there exists a connected double cover $\tilde{N}$ of $N$ such that the lift $\tilde{M}$ of $M$ is a connected orientable embedded hypersurface. Furthermore, $\tilde{M}$ separates $\tilde{N}$ and both components of $\tilde{N}\setminus\tilde{M}$ are diffeomorphic to $N\setminus M$.
\end{proposition}

\subsection{An $H$-radius estimate}

We shall adapt an argument that goes back to Schoen and Yau \cite{SY83} to prove that when taking a lamination limit of closed, embedded minimal surfaces with bounded Morse index the resulting (minimal) lamination only consists of \textsl{proper} leaves. By this we mean that each leaf can be proved to be a closed submanifold in our ambient $(N,g)$. More specifically, we present here an estimate that provides a uniform upper bound for the radius of a stable (intrinsic) disk when the scalar curvature of $g$ is strictly positive, as in the assumptions of our main theorem.

\begin{definition}
Given an open, regular domain $\Omega\subset N$ we shall define its $H-$radius as follows:
\begin{enumerate}
\item{For a simple closed curve $\Gamma$ which bounds an (embedded) open disk in $\Omega$ we set
\[
\textrm{Rad}(\Omega,\Gamma)=\sup\left\{ r: d_g(\Gamma, \partial\Omega)>r, \ \Gamma \ \textrm{does not bound a disk in} \ N_r(\Gamma)\right\}
\] for $N_r(\Gamma)$ the set of points in $N$ within distance $r$ of $\Gamma$;}
\item{at that stage we set
\[
\textrm{Rad}(\Omega)=\sup\left\{ \textrm{Rad}(\Omega, \Gamma) : \Gamma \ \textrm{as above} \right\}.
\]}	
\end{enumerate}			
\end{definition}

Using a suitable weighted arc-length functional and computing its second variation,  Schoen and Yau proved the following.

\begin{proposition}\label{pro:Hrad}\cite{SY83}
Suppose $(N,g)$ is a three-dimensional Riemannian manifold and $\Omega\subset N$ is a bounded region such that the first Dirichlet eigenvalue on $\Omega$ of the operator $-\Delta+\frac{1}{2}R$ is at least $\sigma>0$. Then $\textrm{Rad}(\Omega)\leq \sqrt{\frac{3}{2}}\frac{\pi}{\sqrt{\sigma}}$. 		
\end{proposition}	

We now transpose this method to suit our needs.

\begin{definition}
Let $(S,h)$ be a Riemannian 2-manifold which is complete as a metric space with respect to its own Riemannian distance. If $\Sigma\subset S$ is a simply connected, bounded open domain then one can define its $H-$radius as follows:
\[
\textrm{Rad}(\Sigma)= \sup \left\{ r: \Sigma\setminus N_r(\partial\Sigma)\neq\emptyset    \right\}
\]
\end{definition}

\begin{remark}
	Trivial examples show that the $H-$radius of a domain does not coincide with, nor is equivalent to (half of) its diameter. For instance, if $\Sigma$ is a smoothened rectangle in $\mathbb{R}^2$ then its $H-$radius is comparable with (one half of) the length of its \textsl{shortest} side. However, if (in the setting of the previous definition) $\Sigma$ is a metric ball of radius $\rho$ then it is easily checked that its $H-$radius equals $\rho$. 
\end{remark}	

Getting back to the general setting of our article, let $M\subset N$ be a minimal surface and let $\Sigma\subset M$ be a bounded, regular domain which is diffeomorphic to a disk. If $\Sigma$ is stable, then we notice that the first Dirichlet eigenvalue of the (intrinsic) operator $-\Delta+K$ is at least $\sigma/2$ (because of the standard rearrangement of the Jacobi operator). At this stage, folllowing the argument needed to prove the above proposition one obtains (with purely notational changes) the following assertion, which is of independent interest: 

\begin{proposition}\label{prop:Hradsurf}
Let $(N,g)$ be a three-dimensional manifold whose scalar curvature is bounded below by $\sigma>0$. Let $M$ be a minimally embedded surface in $N$: then any stable disk $\Sigma$ in $M$ has radius bounded from above by $\sqrt{\frac{8}{3}}\frac{\pi}{\sqrt{\sigma}}$. In particular, if $M$ is stable and two-sided, then it is diffeomorphic to $S^2$.		
\end{proposition}	

We have decided to include the complete proof of this proposition, despite the analogies with that of Proposition \ref{pro:Hrad}, both for the convenience of the reader and for the useful simplifications (even at a notational level) that one can obtain in our specific setting.

\begin{proof}
	Let then $\Sigma$ be a bounded, stable simply-connected domain of the minimal surface $M$: then the corresponding Jacobi operator $L$ (with Dirichlet boundary conditions) has negative spectrum\footnote{Let us remark that the domain $\Sigma$, being simply connected by assumption, is necessarily two-sided and hence the use of the scalar stability operator is always legitimate.} and hence, to greater extent, for all test functions $\phi\in C^{\infty}_c(\Sigma)$ one has $\phi J(\phi)\leq 0$ where $J$ is the self-adjoint operator defined by
	\[
	J(\phi)=\Delta\phi+\left(\frac{1}{2}R-K\right)\phi
	\]
	and $K$ denotes the Gauss curvature of $M$. In particular, if we let $\omega$ denote the first (Dirichlet) eigenfunction of the operator J we shall have 
	\begin{equation}\label{eq:stab1}
	J(\omega)=\Delta \omega +\left(\frac{1}{2}R-K\right)\omega\leq 0 
	\end{equation}
	at all points of $\Sigma$. Now, let us pick a number $\rho<\textrm{Rad}(\Sigma)$: then, by the very definition of $\textrm{Rad}(\cdot)$ one can find a point $x\in \Sigma\cap \partial N_{\rho}(\partial\Sigma)$. We can then define, on the class of curves lying in $\Sigma$ and connecting $x$ to a point on the boundary $\partial\Sigma$ the weighted length functional 
	\[
	I(\Gamma)=\int_{\gamma}\omega  ds
	\]
	where $s$ always denotes (here and below) the arc-length parameter of $\gamma$. Said $\gamma$ a curve minimizing such functional, posssibly by replacing it (without renaming) by its terminal segment we can assume it is entirely contained in $N_\rho(\partial\Sigma)$. Now, the second variation of $I(\cdot)$ is defined by a second-order differential operator which reads
	\[
	J_0(\psi)=\frac{d^2\psi}{d s^2}+\omega^{-1}\frac{d\psi}{ds}\frac{d\omega}{ds}+\left(K-\omega^{-1}\Delta \omega+\omega^{-1}\frac{d^2\omega}{ds^2}    \right)\psi
	\]
	which we can rewrite as
	\[
	J_0(\psi)=\frac{d^2\psi}{d s^2}+\omega^{-1}\frac{d\psi}{ds}\frac{d\omega}{ds}+\left(\frac{1}{2}R-\omega^{-1}J \omega+\omega^{-1}\frac{d^2\omega}{ds^2}    \right)\psi
	\]
	and hence, if $h$ is the first eigenfunction of $J_0$ (with Dirichlet boundary conditions on $[0,l]$ for $l=\textrm{length}(\gamma)$) using \eqref{eq:stab1} and the uniform bound $R\geq \sigma$ we derive
	\[
	h^{-1}h''+\omega^{-1}\omega''+h^{-1}\omega^{-1}h'\omega'+\sigma/2 \leq 0.
	\]
	Given $\phi\in C^{\infty}_c[0,l]$, integration by parts allows to write such functional inequality in the form
	\[
	\int_0^{l}\left\{\frac{1}{2}\left(h^{-2}(h')^2+\omega^{-2}(\omega')^2 \right)\phi^2+\frac{1}{2}\left(\frac{d}{ds}\log \omega h  \right)^2\phi^2+\frac{\sigma}{2}\phi^2  \right\}\,ds
	\leq 2\int_{0}^{l}\phi \phi' \frac{d}{ds}\log\omega h \,ds.
	\]
	At this stage, let us exploit the algebraic inequality
	\[
2\left|\phi\phi'\frac{d}{ds}\log \omega h \right|\leq \frac{1}{2}\left(\left(\omega^{-1}\omega'\right)^2+\left(h^{-1}h'\right)^2\right)\phi^2+\frac{1}{2}\left(\frac{d}{ds}\log \omega h \right)^2\phi^2 +\frac{4}{3}(\phi')^2 
	\]
	to obtain the final inequality
	\[
	\frac{\sigma}{2}\int_{0}^{l}\phi^2\,ds\leq \frac{4}{3}\int_{0}^{l}(\phi')^2\,ds
	\]
	which ensures that the operator $-\frac{d^2}{ds^2}-\frac{3}{8}\sigma$ has non-negative first eigenvalue. On the other hand, its Dirichlet spectrum can be computed explicitly and is given by $k^2\pi^2/l^2-3\sigma/8$ for $k=1,2,3,\ldots$ and hence we conclude that necessarily
	\[
	l\leq \sqrt{\frac{8}{3}}\frac{\pi}{\sqrt{\sigma}}
	\]
	 and so, since $l\geq\rho$ and $\rho$ is any number smaller than $\textrm{Rad}(\Sigma)$ this completes the proof.
	 For what concerns the last statement: by virtue of \cite{FCS80} we know that $M$ is conformally equivalent to their the standard sphere $\mathbb{C}\mathbb{P}^{1}$ or the plane $\mathbb{C}$, however the latter alternative cannot possibly occur because of the upper bound on the radius of stable disks and hence $M$ is diffeomorphic to $S^2$. 
\end{proof}

\section{Proof of Theorem \ref{thm:main}}\label{sec:proof}

This section shall almost entirely be devoted to the proof of our main result, with the expection of its last part which contains some concluding remarks concerning the sharpness of Theorem \ref{thm:main}.

\begin{proof}
	
	Let us preliminarily observe that the assumption that $N$ is orientable automatically implies that an embedded surface $M\subset N$ is two-sided if and only if it is orientable.

For the sake of a contradiction, let us then assume the existence of a sequence $\left\{M_i\right\}\subset \mathcal{I}(I)$ (for a fixed integer $I\geq 0$) of pairwise distinct elements. One can then use the index bound as in the recent paper of the B. Sharp \cite{Sha15} to show that given any $M\in \mathcal{I}(I)$ for any small radius $r>0$ there exist at most $I$ balls $B_{r}(x_1),\ldots, B_r(x_I)$ such that
\[
B_r(x)\subset N\setminus \cup_{k=1}^{I}B_r(x_k) \ \Longrightarrow \ M \ \textrm{is stable in} \ B_r(x)
\]
and thus the basic stability estimates by Schoen \cite{Sch83} ensure the existence of a constant $C$ (independent of $M$) such that uniform curvature estimates hold on $B_{r/2}(x)$, namely
\[
\sup_{y\in B_{r/2}(x)}|A(y)|\leq \frac{C}{r}
\]
where $|A(y)|$ denotes the length of the second fundamental form of $M$ at the point $y$.

Considering balls of smaller and smaller radii, a standard diagonal argument leads to the conclusion that (thanks to the lamination compactness statement, Theorem \ref{thm:lamcomp}) a subsequence of $\left\{M_i\right\}$ which we shall not rename converges to a smooth minimal lamination $\mathcal{L}$ of $N\setminus\mathcal{Y}$ away from an exceptional set $\mathcal{Y}=\left\{y_1,\ldots, y_P\right\}$ made of $P\leq I$ points. Furthermore, let us emphasize that (due to the leafwise smooth convergence) we can derive curvature estimates for $\mathcal{L}$ and namely we have that for all $y\in\cup_{\Lambda\in\mathcal{L}}\Lambda$
\[
|A_{\mathcal{L}}(y)|\inf_{y_k\in\mathcal{Y}}d_g(y,y_k)\leq C
\] 
where of course $d_g(\cdot,\cdot)$ denotes the Riemannian distance in $(N,g)$.

Said $\Lambda$ a leaf of $\mathcal{L}$, we recall from Remark \ref{rem:geomHarnack} that the following dichotomy holds: \textsl{either} $\Lambda$ is non-accumulating (which implies that for any $p\in\Lambda$ there is a small cylindrical neighborhood of $p$ in $N$ whose intersection with the support of the lamination $\mathcal{L}$ consists of only one connected component, namely a disk of $\Lambda$) \textsl{or} $\Lambda$ is an accumulating leaf (by which we mean that for any $r>0$ small enough there exists a sequence of connected components $\left\{\Lambda_i\right\}$ of $\cup_{\Lambda'\in\mathcal{L}}\Lambda' \cap B_r(p)$ so that $\Lambda_i\to\Lambda$ smoothly with multiplicity one). 


Now, the cases when i) either $\mathcal{L}$ contains an accumulating leaf $\Lambda$ or ii) all the leafs of $\mathcal{L}$ are non-accumulating, but there is at least one leaf $\Lambda$ such that $M_i\to \Lambda$ smoothly and graphically with multiplicity greater than one\footnote{By virtue of the lamination convergence, for any non-concentration point $p\in\Lambda$ there exists a small cylindrical neighbourhood $\Gamma$ (centered at $p$) such that $M_i\cap \Gamma$ consists of finitely many (say $\chi(i)$) connected components, which are in fact graphs over the base $\Gamma\cap\Lambda$ with uniformly bounded slope for $i\geq i_0=i_0(\Gamma)$: thus, possibly extracting a further subsequence, we can assume that either $\chi(i)=1$ for all $i\geq i_0$ or $\chi(i)\geq 2 $ for all $i\geq i_0$. In the former case, we shall say that the convergence happens with multiplicity one, while in the latter we shall say that the convergence happens with mulitplicity greater than one.} can be handled in a rather similar fashion. 

Let us first assume that $\Lambda$ is two-sided so that (once a choice of a unit normal field $\nu$ is made) for any regular, connected open set $\Omega\subset\subset \Lambda$ we can fix $\eta>0$ so that the intersection of $M_i$ with the tubular neighborhood $U$ of size $\eta$ over $\Omega$ will eventually consist of $\chi(i)$ graphs\footnote{With possibly (in the first case) some additional connected components which are not globally graphical over $\Omega$ and lie outside of the slab bounded below by the graph of $u^{1}_{i}$ and above by $u^{\chi(i)}_i$.} with smooth defining functions $u^{1}_i<u^2_{i}<\ldots< u^{\chi(i)}_{i}$
where in the former case $\chi(i)\to\infty$ as $i\to\infty$, while in the latter $\chi(i)$ is uniformly bounded below by two. 
For a given, fixed, point $p\in\Omega$, let $n(i)$ be the positive integer denoting the couple of graphs whose heights above $p$ are minimal, more precisely
\[
n(i)=\max\left\{ n \ | \ 1\leq n<\chi(i),  n\in\arg\min_{1\leq m<\chi(i)} |u^{m}_{i}(p)|+|u^{m+1}_i(p)|\right\}
\]
where taking the maximum (in the class of integers for which such quantity achieves the least value\footnote{Notice that this class has either cardinality one or two.}) is only needed to avoid ambiguities for symmetric configurations.
Then, if we set   
\[
w^{(\Omega)}_{i}=u^{n(i)+1}_{i}-u^{n(i)}_{i}, \ \ u^{(\Omega)}_i=w_i/w_i(p)
\]
 the function $u^{(\Omega)}_i:\Omega\to\mathbb{R}$ is well-defined (at least for $i\geq i_0, i_0=i_0(\Omega)$), smooth and positive.
 In both cases obviously $w_i(p)>0$ by embeddedness of each surface $M_i$ and in fact $w_i>0$ at all points).
 Exploiting the minimality of $\Lambda$ and of $M_i$ one can then show, as in \cite{Sim87} (or e.g. page 11 in \cite{Sha15}), that (possibly extracting a further subsequence) linear Harnack and elliptic estimates ensure $u_i\to u$ with $u(p)=1$ and $Lu=0$ for $L$ the Jacobi operator of $\Lambda$, namely $u$ is a positive solution of the Jacobi operator. This puts us in the condition of applying Barta's criterion (cmp. for instance Lemma 1.36 in \cite{CM11}) to conclude that $\Lambda$ is stable. As a result, applying the removable singularity Theorem \ref{thm:remsing} at most $I$ times\footnote{Notice that, by construction, $\mathcal{L}$ is a minimal lamination of $N\setminus\mathcal{Y}$ and thus $\Lambda$ can also be regarded as a smooth, minimal lamination in small \textsl{punctured} balls centered at $y_k$ for $1\leq k\leq I$.}, we can show that $\Lambda$ extends to a complete, embedded minimal surface $\Lambda_{\mathcal{Y}}$ in $(N,g)$. 

That being achieved, we can prove that the function $u$ in question extends across the (removed) punctures $\mathcal{Y}$ of $\Lambda_{\mathcal{Y}}$, which follows via standard local theory of elliptic PDEs by proving uniform $L^{\infty}$ bounds for the functions $u_i$: specifically we shall observe that the supremum of $u_i$ in a suitable (finite) solid cylinder $\Gamma_{\varepsilon}(y_k)$ is controlled by a (universal constant times) the supremum of the same function on the surface $\partial \Gamma_\varepsilon$ of the same cylinder for some $\varepsilon>0$ small enough. This is accomplished in two steps, by first proving that $M_i\to \mathcal{L}$ in Hausdorff distance (which is a trivial argument by contradiction) and then using the local \textsl{minimal foliation trick}
due to B. White \cite{Whi87} and the maximum principle as is done in Appendix A of \cite{CM00}. So, in the end, we have constructed a positive solution $\overline{u}$ of the Jacobi equation for the complete minimal surface $\Lambda_{\mathcal{Y}}$. 

Lastly, let us observe that the positive lower bound on the ambient scalar curvature implies, thanks to well-known results in \cite{FCS80} that $\Lambda$, being stable, is diffeomorphic (in fact conformally equivalent to) either a sphere of a plane. The latter case is immediately ruled out by means of our $H-$radius estimates, Proposition \ref{prop:Hradsurf} and so $\Lambda$ must be a stable sphere inside $(N,g)$.

By virtue of the previous steps we already know that such minimal sphere has a non-trivial Jacobi vector field, which in turn contradicts the bumpyness assumption on the Riemannian metric $g$. 

Let us now discuss the necessary modifications to such argument in case the leaf $\Lambda$ in question is instead one-sided. First of all, we notice that given $y_k\in\mathcal{Y}$ (the usual exceptional set for the convergence) there exists $r>0$ small enough so that $\Lambda\cap B_r(y_k)\setminus\left\{y_k\right\}$ is orientable with trivial normal bundle (cmp. Proposition 4.2 in \cite{LZ15}) hence (possibly taking an even smaller $r$) the same conclusion will hold for all connected components of $\mathcal{L}$ in $B_r(p)$ and, in addition, $\Lambda\cap B_r(y_k)\setminus\left\{y_k\right\}$  is stable, by \textsl{local} construction of a positive solution of the Jacobi equation. Hence we can apply Theorem \ref{thm:remsing} to ensure that $\Lambda$ extends in fact to a smooth minimal leaf (which we shall rename, as usual, $\Lambda_{\mathcal{Y}}$).

Then we can use Proposition \ref{pro:xz} and construct an orientable double cover $\pi:\tilde{N}\to N$ such that the pre-image $\tilde{\Lambda}_{\mathcal{Y}}=\pi^{-1}(\Lambda_{\mathcal{Y}})$ is two-sided (and separating). Notice that this does not rely on the compactness of $\Lambda_{\mathcal{Y}}$ (namely this construction does not need $\Lambda_{\mathcal{Y}}$ to be a closed surface). In particular, the lifted leaf $\tilde{\Lambda}_{\mathcal{Y}}$ has a well-defined unit normal $\nu$ which we can use to unambigously define, near any point $x\in \tilde{\Lambda}$ the graphical sheets of $\pi^{-1}(M_i)$ (locally, in a tubular neighborhood of $\tilde{\Lambda}_{\mathcal{Y}}$) and hence, arguing as in the simply-connected case, we can then produce a positive solution of the Jacobi equation for $\tilde{\Lambda}_{\mathcal{Y}}$. Thus $\tilde{\Lambda}_{\mathcal{Y}}$ has to be stable, hence by the $H-$radius estimate provided by Proposition \ref{prop:Hradsurf} it has to be compact and so necessarily a sphere since we have a positive bound on the scalar curvature of $(N,g)$. It follows that $\Lambda_{\mathcal{Y}}$ is diffeomorphic to $\mathbb{R}\mathbb{P}^2$.

\

This concludes the proof in the cases listed above, and leaves us with the treatment of the case when $\left\{M_i\right\}$ subconverges to any leaf of $\mathcal{L}$ with multiplicity one. 

\

Let $\Lambda$ be a leaf of $\mathcal{L}$ (in fact, in this case we will see below that $\mathcal{L}$ must consist of only one leaf, namely such $\Lambda$). 
Once again, let us start by considering the simpler case when $\Lambda$ is two-sided.

In this case we can articulate the argument as follows: 
\begin{enumerate}
\item[a)]{there exists $\varepsilon$ such that for each $y_k, \ k=1,\ldots, P$ in the set of exceptional points $\mathcal{Y}$: if $y_k\in \mathcal{Y}\cap(\overline{\Lambda}\setminus\Lambda)$ then such leaf $\Lambda$ is stable in the punctured ball $B_{\varepsilon}(y_k)\setminus\left\{y_k\right\}$ (this follows, by contradiction, as in \cite{Sha15} for otherwise one could exploit the strong single-sheeted convergence to violate the index bound for $M_i, \ \textrm{for} \ i \ \textrm{large}$, by constructing a sequence of area-decreasing vector fields supported in disjoint annuli centered at $y_k$);}
\item[b)]{by Theorem \ref{thm:remsing} $\Lambda$ as in a) has a removable singularity at $y_k$, so that $\Lambda$ extends to a smooth, complete embedded minimal surface $\Lambda_{\mathcal{Y}}$ in $(N,g)$;}
\item[c)]{we claim that all leafs of $\mathcal{L}$ have to be closed surfaces. If, on the contrary, $\Pi$ were a (possibly desingularized) complete, non properly embedded leaf of $\mathcal{L}$ diffeomorphic to a punctured Riemann surface, then by compactness of $N$ one could find $p\in N, p_0\in\Pi$ and a sequence $\left\{p_j \right\}\subset \Pi$ such that 
	\begin{enumerate}
	\item[i)]{said $d_{\Pi}(\cdot,\cdot\cdot)$ the intrinsic distance induced on $\Pi$ one has $d_{\Pi}(p_0,p_j)\nearrow+\infty$;}	
	\item[ii)]{in $(N,g)$ the sequence $\left\{p_j\right\}$ converges to $p$.}	
	\end{enumerate}	
	(Notice that since $\mathcal{L}$ is closed (by definition of lamination) necessarily $p$ belongs to a leaf, say $p\in\Lambda^{\omega}_{\mathcal{Y}}$.) 
By assumption, the lamination $\mathcal{L}$ is postulated in this case not to have accumulating leaves, and so in particular $\Lambda^{\omega}_{\mathcal{Y}}$ is a non-accumulating leaf, which implies that one can find $r>0$ such that $\Lambda^{\omega}_{\mathcal{Y}}\cap B_r(p)$ consists of only one connected component (namely a disk containing $p$). Hence,  necessarily the point $p_j$ belongs to such component for all $j$ large enough. Thus $d_{\Pi}(\cdot,\cdot\cdot)$ must be uniformly bounded, contradiction;}
\item[d)]{by b) and c) each leaf of $\mathcal{L}$ must be diffeomorphic to a closed Riemann surface, thus the Hausdorff convergence of $M_{i}\to\mathcal{L}$ plus a standard connectedness argument ensure that $\mathcal{L}$ consists of only one closed leaf $\Lambda_{\mathcal{Y}}$;}
\item[e)]{arguing as in the first part of the proof we can construct a non-trivial function $u:\Lambda_{\mathcal{Y}}\to\mathbb{R}$ which solves $Lu=0$ for $L$ the Jacobi operator of $\Lambda_{\mathcal{Y}}$ (even though such function does not need to be positive). The fact that the function in question does not vanish identically follows from the Hausdorff convergence (gained in part d)), the existence of a local foliation made of minimal leaves and the maximum principle as in the very last paragraph of the proof of Theorem 1.1 in \cite{CM00}: if we let $\left\{\Omega_i\right\}$ be an exhaustion of $\Lambda$ and $u_i:\Omega_i\to\mathbb{R}$ be the defining function of $M_i$ (with $\Omega_i$ chosen so that this is well-defined) then one can prove that for every fixed $\delta>0$ and $p\in\Lambda_{\mathcal{Y}}\setminus\Lambda$ one has a bound of the form
		\[
		\sup_{x\in B_{\delta/2}(p)\cap\Omega_i} |u_i(x)|\leq C \sup_{x\in\partial B_{\delta}(p)}|u_i(x)|
		\]  
which ensures that, set $h_i(x):=u_i(x)/|u_i|_{L^2(\Omega_i)}$ then it is not possible that $h_i\to 0$ on $\partial B_{\delta}(p)$ since this would necessarily violate the normalization condition we have imposed, namely $|h_i|_{L^2(\Omega_i)}=1$. Lastly, the very same bounds allow to smoothly extend the function $u$ from $\Lambda$ to $\Lambda_{\mathcal{Y}}$;}
\item[f)]{we reach a contradiction because of the bumpyness assumption on the metric $g$.}
\end{enumerate}	Let us conclude our proof by handling this multiplicity one case when $\Lambda$ is one-sided. Arguing locally, we can remove the singularities (if any) by exploiting the index bound on $M_i$ by virtue of the multiplicity one convergence.
	
	  Let us lift our objects to the double cover $\pi:\tilde{N}\to N$ as above. Due to the multiplicity one convergence (which is related to the fact that the local scenario in $(N,g)$ is isometric to that in $(\tilde{N},\tilde{g})$), for $i$ large enough we can express $\tilde{M}_i$ as a normal graph over $\tilde{\Lambda}_{\mathcal{Y}}$ by means of the exponential map, say $u_i:\tilde{\Lambda}_{\mathcal{Y}}\to\mathbb{R}$ be the corresponding defining function. The relation $\tau_{\ast} (\nu)=-\nu$ implies that $u_i(x)=-u_i(\tau(x))$ for every $x\in\tilde{\Lambda}_{\mathcal{Y}}$. 
	  We can then find a smooth function $\tilde{u}:\tilde{\Lambda}_{\mathcal{Y}}\to\mathbb{R}$ such that $\tilde{u}(\tau(x))=-\tilde{u}(x)$ solving $\tilde{L}\tilde{u}=0$. It follows that the vector field $\tilde{v}=\tilde{u}\nu$ satisfies $\tau_{\ast}(\tilde{v})=\tilde{v}$ and hence determies a well-defined Jacobi field for $\Lambda_{\mathcal{Y}}\subset N$, thus leading to a contradiction.

Thereby, the proof is complete.
	
\end{proof}	

\begin{remark}\label{rem:counter}
Given the argument above, it seems appropriate to add a few comments about the role of the assumption that the manifold $(N,g)$ does not contain any minimally embedded $\mathbb{R}\mathbb{P}^2$ for the generic finiteness result to be true. 
Throughout this remark, let us identify $S^3$ (as a smooth manifold) with the unit sphere in $\mathbb{R}^4$.
Let us start by saying that it is easy to construct an \textsl{even} Riemannian metric $g$ on $S^3$ of positive scalar curvature such that $(S^3,g)$ contains a sequence of \textsl{strictly stable} minimal spheres $\left\{M_i\right\}$ such that:
\begin{enumerate}
\item{the areas of these surfaces are (eventually) strictly increasing in $i$
	\[
\mathcal{H}^2(M_i)<\mathcal{H}^2(M_{i+1})<\mathcal{H}^2(M_{i+2})<\ldots
\]
for all $i\geq i_0$ a given large integer;}
\item{the surfaces in question converge strongly (with multiplicity one and from one side) to a totally geodesic sphere $M_{\infty}$.}
\end{enumerate}
The metric $g$ is (the extension of) a warped-product metric on $S^{3}\setminus\left\{NP,SP\right\}$ where of course the north pole $NP$ has cordinates $(0,0,0,1)$ and the south pole $SP$ has coordinates $(0,0,0,-1)$. 
More precisely, given any $\varepsilon>0$ we can pick $\delta>0$ in a way that the metric on $S^2\times (-1,1)$ given by $g(x,t)=f(t)^2g_{S^2}(x)+dt\otimes dt$ for $f\in C^{\infty}((-1,1);\mathbb{R})$ with $f(t)=f(-t)$ and
\[
f(t)=\begin{cases}
1+e^{-1/t}\sin\left(\frac{1}{t}\right) & \textrm{if} \ t\in(0,\delta) \\
\sin(t) & \textrm{if} \ t\in(2\delta,1)
\end{cases}
\]	
has scalar curvature
\[
R(g_{\varepsilon})=-4\frac{f''}{f}+2\frac{1-(f')^2}{f^2}
\]
 bounded from below by $2-\varepsilon$. The minimal surfaces of the form $M:=\left\{t=c\right\}$ for $c\in(-1,1)$ correspond to the critical points of the warping factor $f$: in particular if we set $M_{i}:=\left\{t=\left(\frac{\pi}{4}+i\pi \right)^{-1} \right\}$ (for $i$ large enough) we obtain a sequence of embedded minimal spheres which are strictly stable provided $i$ is odd. Furthermore
\[
\mathcal{H}^2(M_i)=4\pi \left(1-\frac{e^{-(\frac{\pi}{4}+i\pi)}}{\sqrt{2}} \right)^2
\]
which is patently increasing in $i$. 
Each of those metrics is of course not bumpy, but we expect it should be possible to modify them in the sole spherical cap $S^2\times (2\delta, 1)$ (and, symmetrically, on $S^2\times (-1,-1+2\delta)$) to obtain new even metrics $\tilde{g}_{\varepsilon}$ such that the only minimal surface in $(S^3,\tilde{g}_{\varepsilon})$ having a Jacobi field is the equatorial (totally geodesic) $S^2$ given by $M_{\infty}:=\left\{t=0 \right\}$ with such Jacobi field being of the form $v=\nu$ for $\nu$ a choice of the unit normal to $M_{\infty}$. Now, if that were the case, we could of course project this whole picture to $\mathbb{R}\mathbb{P}^3$ thereby obtaining a closed 3-manifold (endowed with a bumpy Riemannian metric whose scalar curvature bounded from above by $2-\varepsilon$) that contains infinitely many strictly stable minimal spheres. 
\end{remark}	

\

\end{document}